\documentclass[11pt]{amsart}
\usepackage{amsmath,amssymb}
\newtheorem{theorem}{Theorem}[section]
\newtheorem{proposition}[theorem]{Proposition}
\newtheorem{corollary}[theorem]{Corollary}
\newtheorem{lemma}[theorem]{Lemma}
\theoremstyle{definition}
\newtheorem{definition}[theorem]{Definition}

\begin{document}

\title[Mean curvature flow with surgery]{A monotonicity formula for mean curvature flow with surgery} 
\author{Simon Brendle}
\address{Department of Mathematics \\ Stanford University \\ Stanford, CA 94305}
\thanks{The author was supported in part by the National Science Foundation under grant DMS-1201924.}
\begin{abstract}
We prove a monotonicity formula for mean curvature flow with surgery. This formula differs from Huisken's monotonicity formula by an extra term involving the mean curvature. As a consequence, we show that a surgically modified flow which is sufficiently close to a smooth flow in the sense of geometric measure theory is, in fact, free of surgeries. This result has applications to the longtime behavior of mean curvature flow with surgery in Riemannian three-manifolds.
\end{abstract}
\maketitle 

\section{Introduction}

One of the main tools in the study of mean curvature flow is Huisken's monotonicity formula (cf. \cite{Huisken}). In the special case of smooth solution of mean curvature flow in $\mathbb{R}^3$, the monotonicity formula implies that the Gaussian integral 
\[\int_{M_t} \frac{1}{4\pi (t_0-t)} \, e^{-\frac{|x-p|^2}{4(t_0-t)}}\] 
is monotone decreasing in $t$, as long as $t_0-t>0$. This mononotonicity property plays a crucial role in the analysis of singularities; see e.g. \cite{Colding-Minicozzi}, \cite{Haslhofer-Kleiner1}, \cite{White1}, \cite{White2}, \cite{White3}.

Our goal in this paper is to adapt the monotonicity formula to surgically modified flows. Motivated by earlier work of Hamilton \cite{Hamilton1},\cite{Hamilton2} and Perelman \cite{Perelman1},\cite{Perelman2} on the Ricci flow, Huisken and Sinestrari \cite{Huisken-Sinestrari} introduced a notion of mean curvature flow with surgery for two-convex hypersurfaces in $\mathbb{R}^{n+1}$, where $n \geq 3$. In a joint work with Gerhard Huisken \cite{Brendle-Huisken1}, we extended this construction to the case $n=2$. As a result, we obtained a notion of mean curvature flow with surgery for embedded, mean convex surfaces in $\mathbb{R}^3$. One of the key ingredients in the proof is a sharp estimate for the inscribed radius established earlier in \cite{Brendle1}. An alternative construction was given by Haslhofer and Kleiner \cite{Haslhofer-Kleiner2}. We note that the surgery construction in \cite{Brendle-Huisken1} can be extended to flows of embedded, mean convex surfaces in three-manifolds; see \cite{Brendle2} and \cite{Brendle-Huisken2} for details.

Throughout this paper, we will focus on mean curvature flow with surgery for embedded, mean convex surfaces in three-manifolds. It turns out that there is a monotonicity property for surgically modified flows which is similar to the one established by Huisken in the smooth case; however, we need to include an extra term in order to ensure that the monotonicity holds across surgery times. More precisely, if $M_t$ is a mean curvature flow with surgery in $\mathbb{R}^3$, we show that the quantity 
\[\int_{M_t} \frac{1}{4\pi (t_0-t)} \, e^{-\frac{|x-p|^2}{4(t_0-t)}-\frac{H}{200 \, H_1}}\] 
is monotone decreasing in $t$, as long as $t_0-t \geq \frac{5}{9} \, H_1^{-2}$. Here, $H_1$ is a positive constant which represents the so-called surgery scale; in other words, each neck on which we perform surgery has radius between $\frac{1}{2H_1}$ and $\frac{1}{H_1}$. A similar monotonicity property holds for mean curvature flow with surgery in a Riemannian three-manifold.

In Section \ref{auxiliary.results}, we establish a number of auxiliary results. These results will be used in Section \ref{monotonicity.for.mcf.with.surgery} to deduce a monotonicity formula for mean curvature flow with surgery in $\mathbb{R}^3$. In Section \ref{Riem.manifolds}, we extend the monotonicity formula to solutions of mean curvature flow with surgery in Riemannian three-manifolds. 

\section{Behavior of a Gaussian integral under a single surgery}

\label{auxiliary.results}

\begin{lemma}
\label{technical.ingredient.1}
There exists a real number $\beta>0$ with the following significance. Suppose that $\Gamma$ is a curve in $\mathbb{R}^2$ which is $\beta$-close to the unit circle in the $C^1$-norm. Moreover, suppose that $\psi$ is a real-valued function defined on $\Gamma$ satisfying $\sup_\Gamma |\psi-1| < \beta$. Then 
\[\int_\Gamma \psi \, e^{-\frac{|x|^2}{4\tau}+\langle q,x \rangle} \, \Big ( \frac{10}{11} - \frac{|x|^2}{2\tau} + \langle q,x \rangle \Big ) \geq 0\] 
for all $\tau \geq \frac{5}{9}$ and all $q \in \mathbb{R}^2$.
\end{lemma}

\textbf{Proof.} 
Let us fix a large constant $Q$ with the property that 
\[e^{\frac{|q|}{2}} \int_{\Gamma \cap \{\langle q,x \rangle \geq \frac{|q|}{2}\}} \psi \geq 16 \int_{\Gamma \cap \{\langle q,x \rangle \leq 0\}} \psi\] 
for all points $q \in \mathbb{R}^2$ satisfying $|q| \geq Q$. This implies 
\[e^{\frac{|q|}{2}} \int_{\Gamma \cap \{\langle q,x \rangle \geq \frac{|q|}{2}\}} \psi \, e^{-\frac{|x|^2}{4\tau}} \geq 4 \int_{\Gamma \cap \{\langle q,x \rangle \leq 0\}} \psi \, e^{-\frac{|x|^2}{4\tau}}\] 
for all $\tau \geq \frac{5}{9}$ and all $q \in \mathbb{R}^2$ satisfying $|q| \geq Q$. From this, we deduce that  
\begin{align*} 
\int_{\Gamma \cap \{\langle q,x \rangle \geq \frac{|q|}{2}\}} \psi \, e^{-\frac{|x|^2}{4\tau}+\langle q,x \rangle} \, \langle q,x \rangle 
&\geq \frac{|q|}{2} \, e^{\frac{|q|}{2}} \int_{\Gamma \cap \{\langle q,x \rangle \geq \frac{|q|}{2}\}} \psi \, e^{-\frac{|x|^2}{4\tau}} \\ 
&\geq 2 \, |q| \int_{\Gamma \cap \{\langle q,x \rangle \leq 0\}} \psi \, e^{-\frac{|x|^2}{4\tau}} \\ 
&\geq -\int_{\Gamma \cap \{\langle q,x \rangle \leq 0\}} \psi \, e^{-\frac{|x|^2}{4\tau}+\langle q,x \rangle} \, \langle q,x \rangle 
\end{align*}
for all $\tau \geq \frac{5}{9}$ and all $q \in \mathbb{R}^2$ satisfying $|q| \geq Q$. Therefore, we obtain 
\[\int_\Gamma \psi \, e^{-\frac{|x|^2}{4\tau}+\langle q,x \rangle} \, \langle q,x \rangle \geq 0,\] 
hence 
\[\int_\Gamma \psi \, e^{-\frac{|x|^2}{4\tau}+\langle q,x \rangle} \, \Big ( \frac{10}{11} - \frac{|x|^2}{2\tau} + \langle q,x \rangle \Big ) \geq 0\] 
for all $\tau \geq \frac{5}{9}$ and all $q \in \mathbb{R}^2$ satisfying $|q| \geq Q$. On the other hand, by choosing $\beta>0$ sufficiently small, we can arrange that 
\[\int_\Gamma \psi \, e^{-\frac{|x|^2}{4\tau}+\langle q,x \rangle} \, \Big ( \frac{10}{11} - \frac{|x|^2}{2\tau} + \langle q,x \rangle \Big ) \geq 0\] 
for all $\tau \geq \frac{5}{9}$ and all $q \in \mathbb{R}^2$ satisfying $|q| \leq Q$. This proves the assertion. \\

\begin{lemma}
\label{technical.ingredient.2}
Let $\beta>0$ be chosen as in Lemma \ref{technical.ingredient.1}. Suppose that $\Gamma$ is a curve in $\mathbb{R}^2$ which is $\beta$-close to the unit circle in the $C^1$-norm. Moreover, suppose that $\psi$ is a real-valued function defined on $\Gamma$ satisfying $\sup_\Gamma |\psi-1| < \beta$. Then 
\[\rho^{\frac{10}{11}} \int_\Gamma \psi \, e^{-\frac{\rho^2 \, |x|^2}{4\tau}+\rho \, \langle q,x \rangle} \leq 
\int_\Gamma \psi \, e^{-\frac{|x|^2}{4\tau}+\langle q,x \rangle}\] 
for all $\rho \in (0,1)$, all $\tau \geq \frac{5}{9}$, and all $q \in \mathbb{R}^2$.
\end{lemma}

\textbf{Proof.} 
Let us fix a real number $\tau \geq \frac{5}{9}$ and a point $q \in \mathbb{R}^2$. Moreover, let $\rho \in (0,1)$. Applying Lemma \ref{technical.ingredient.1} with $\tilde{\tau} = \frac{\tau}{\rho^2} \geq \frac{5}{9}$ and $\tilde{q} = \rho \, q$ gives 
\[\int_\Gamma \psi \, e^{-\frac{\rho^2 \, |x|^2}{4\tau}+\rho \, \langle q,x \rangle} \, \Big ( \frac{10}{11} - \frac{\rho^2 \, |x|^2}{2\tau} + \rho \, \langle q,x \rangle \Big ) \geq 0\] 
for all $\rho \in (0,1)$. Therefore, we obtain 
\begin{align*} 
\rho \, \frac{d}{d\rho} \bigg ( \int_\Gamma \psi \, e^{-\frac{\rho^2 \, |x|^2}{4\tau}+\rho \, \langle q,x \rangle} \bigg ) 
&= \int_\Gamma \psi \, e^{-\frac{\rho^2 \, |x|^2}{4\tau}+\rho \, \langle q,x \rangle} \, \Big ( -\frac{\rho^2 \, |x|^2}{2\tau} + \rho \, \langle q,x \rangle \Big ) \\ 
&\geq -\frac{10}{11} \, \bigg ( \int_\Gamma \psi \, e^{-\frac{\rho^2 \, |x|^2}{4\tau}+\rho \, \langle q,x \rangle} \bigg ) 
\end{align*} 
for all $\rho \in (0,1)$. Consequently, the function 
\[\rho \mapsto \rho^{\frac{10}{11}} \int_\Gamma \psi \, e^{-\frac{\rho^2 \, |x|^2}{4\tau}+\rho \, \langle q,x \rangle}\] 
is monotone increasing for $\rho \in (0,1)$. From this, the assertion follows. \\

In the remainder of this section, we consider an $(\hat{\alpha},\hat{\delta},\varepsilon,L)$-neck $N$ of size $1$ (see \cite{Brendle-Huisken1} for the definition). It is understood that $\varepsilon$ is much smaller than $\hat{\delta}$. By definition, we can find a simple closed, convex curve $\Gamma$ with the property that $\text{\rm dist}_{C^{20}}(N,\Gamma \times [-L,L]) \leq \varepsilon$. Since $\text{\rm dist}_{C^{20}}(N,\Gamma \times [-L,L]) \leq \varepsilon$, we can find a collection of curves $\Gamma_s$ such that 
\[\{(\gamma_s(t),s): s \in [-(L-1),L-1], \, t \in [0,1]\} \subset N\] 
and 
\[\sum_{k+l \leq 20} \Big | \frac{\partial^k}{\partial s^k} \, \frac{\partial^l}{\partial t^l} (\gamma_s(t) - \gamma(t)) \Big | \leq O(\varepsilon).\]
Here, we have used the notation $\Gamma = \{\gamma(t): t \in [0,1]\}$ and $\Gamma_s = \{\gamma_s(t): t \in [0,1]\}$.

As in \cite{Brendle-Huisken1}, we may translate the neck $N$ in space so that the center of mass of $\Gamma$ is at the origin. Using the curve shortening flow, we can construct a homotopy $\tilde{\gamma}_r(t)$, $(r,t) \in [0,1] \times [0,1]$, with the following properties: 
\begin{itemize} 
\item $\tilde{\gamma}_r(t)=\gamma(t)$ for $r \in [0,\frac{1}{4}]$.
\item $\tilde{\gamma}_r(t)=(\cos (2\pi t),\sin(2\pi t))$ for $r \in [\frac{1}{2},1]$.
\item For each $r \in [0,1]$, the curve $\tilde{\Gamma}_r$ is $\frac{1}{1+\hat{\delta}}$-noncollapsed.
\item We have $\sup_{(r,t) \in [0,1] \times [0,1]} |\frac{\partial}{\partial r} \tilde{\gamma}_r(t)| + |\frac{\partial^2}{\partial r \, \partial t} \tilde{\gamma}_r(t)| + |\frac{\partial^2}{\partial r^2} \tilde{\gamma}_r(t)| \leq \omega(\hat{\delta})$, where $\omega(\hat{\delta}) \to 0$ as $\hat{\delta} \to 0$.
\end{itemize} 
Finally, we choose a smooth cutoff function $\chi: \mathbb{R} \to \mathbb{R}$ such that $\chi = 1$ on $(-\infty,1]$ and $\chi=0$ on $[2,\infty)$. We next define a surface $\tilde{F}_\Lambda: [-L,\Lambda] \times [0,1] \to \mathbb{R}^3$ by 
\[\tilde{F}_\Lambda(s,t) = \begin{cases} (\gamma_s(t),s) & \text{\rm for $s \in [-(L-1),0]$} \\ ((1-e^{-\frac{4\Lambda}{s}}) \, \gamma_s(t),s) & \text{\rm for $s \in (0,\Lambda^{\frac{1}{4}}]$} \\ ((1-e^{-\frac{4\Lambda}{s}}) \, (\chi(s/\Lambda^{\frac{1}{4}}) \, \gamma_s(t)+(1-\chi(s/\Lambda^{\frac{1}{4}})) \, \gamma(t)),s) & \text{\rm for $s \in (\Lambda^{\frac{1}{4}},2 \, \Lambda^{\frac{1}{4}}]$} \\ ((1-e^{-\frac{4\Lambda}{s}}) \, \tilde{\gamma}_{s/\Lambda}(t),s) & \text{\rm for $s \in (2 \, \Lambda^{\frac{1}{4}},\Lambda]$.} \end{cases}\] 
It is clear that $\tilde{F}_\Lambda$ is smooth. Moreover, $\tilde{F}_\Lambda$ is axially symmetric for $s \geq \frac{\Lambda}{2}$. As in \cite{Brendle-Huisken1}, we may extend the immersion $\tilde{F}_\Lambda$ by gluing in an axially symmetric cap. To do that, we fix a smooth, convex, even function $\Phi: \mathbb{R} \to \mathbb{R}$ such that $\Phi(z) = |z|$ for $|z| \geq \frac{1}{100}$. We then define 
\[a = 1-e^{-4}+\frac{1}{3} \, (1-e^{-4})^2 \, \Lambda^{-\frac{1}{4}}\]
and 
\begin{align*} 
v_\Lambda(s) 
&= 1-e^{-\frac{4\Lambda}{s}} + a \, \sqrt{\frac{\Lambda+2\,\Lambda^{\frac{1}{4}}-s}{a+\Lambda+2\,\Lambda^{\frac{1}{4}}-s}} \\ 
&- \Lambda^{-\frac{1}{4}} \, \Phi \bigg ( \Lambda^{\frac{1}{4}} \, \Big ( 1-e^{-\frac{4\Lambda}{s}} - a \, \sqrt{\frac{\Lambda+2\,\Lambda^{\frac{1}{4}}-s}{a+\Lambda+2\,\Lambda^{\frac{1}{4}}-s}} \Big ) \bigg ) 
\end{align*} 
for $s \in [\Lambda,\Lambda+\Lambda^{\frac{1}{4}}]$. Moreover, we put 
\[v_\Lambda(s) = 2a \, \sqrt{\frac{\Lambda+2\,\Lambda^{\frac{1}{4}}-s}{a+\Lambda+2\,\Lambda^{\frac{1}{4}}-s}}\] 
for $s \in [\Lambda+\Lambda^{\frac{1}{4}},\Lambda+2\Lambda^{\frac{1}{4}}]$. The axially symmetric cap is chosen so that its cross section at height $s \in [\Lambda,\Lambda+2\Lambda^{\frac{1}{4}}]$ is a circle of radius $\frac{1}{2} \, v_\Lambda(s)$. The resulting surface will be denoted by $\tilde{N}$. 

We first consider the region $s \in (0,\Lambda^{\frac{1}{4}}]$.

\begin{lemma}
\label{x_3.small}
There exist positive real numbers $\delta_*$ and $\Lambda_*$ with the following significance. If $\hat{\delta} < \delta_*$ and $\Lambda > \Lambda_*$, then we have 
\begin{align*} 
&\int_{\tilde{N} \cap \{x_3=s\}} e^{-\frac{x_1^2+x_2^2}{4\tau}+q_1x_1+q_2x_2-r_0 \tilde{H}} \, \frac{1}{|\nabla^{\tilde{N}} x_3|} \\ 
&\leq \int_{N \cap \{x_3=s\}} e^{-\frac{x_1^2+x_2^2}{4\tau}+q_1x_1+q_2x_2-r_0 H} \, \frac{1}{|\nabla^N x_3|} 
\end{align*}
for all $s \in (0,\Lambda^{\frac{1}{4}}]$, all $r_0 \in [\frac{1}{1000},1]$, all $\tau \geq \frac{5}{9}$, and all $(q_1,q_2) \in \mathbb{R}^2$.
\end{lemma}

\textbf{Proof.} 
For each $s \in (0,\Lambda^{\frac{1}{4}}]$, the cross section $\tilde{N} \cap \{x_3=s\}$ is obtained by dilating the cross section $N \cap \{x_3=s\}$ by the factor $\rho(s) := 1-e^{-\frac{4\Lambda}{s}}$. Moreover, we have the pointwise inequalities
\[\frac{|\nabla^N x_3|}{|\nabla^{\tilde{N}} x_3|} \leq 1+C \, \frac{\Lambda}{x_3^2} \, e^{-\frac{4\Lambda}{x_3}}\] 
and 
\[\tilde{H}-H \geq c \, \frac{\Lambda^2}{x_3^4} \, e^{-\frac{4\Lambda}{x_3}}\] 
for $x_3 \in (0,\Lambda^{\frac{1}{4}}]$. Hence, if we choose $\Lambda$ sufficiently large, then we have 
\[e^{-r_0 \tilde{H}} \, \frac{1}{|\nabla^{\tilde{N}} x_3|} \leq e^{-r_0 H} \, \frac{1}{|\nabla^N x_3|}\] 
for $x_3 \in (0,\Lambda^{\frac{1}{4}}]$. From this, we deduce that 
\begin{align*} 
&\int_{\tilde{N} \cap \{x_3=s\}} e^{-\frac{x_1^2+x_2^2}{4\tau}+q_1x_1+q_2x_2-r_0 \tilde{H}} \, \frac{1}{|\nabla^{\tilde{N}} x_3|} \\ 
&\leq \rho(s) \int_{N \cap \{x_3=s\}} e^{-\frac{\rho(s)^2 \, (x_1^2+x_2^2)}{4\tau}+\rho(s) \, (q_1x_1+q_2x_2)-r_0 H} \, \frac{1}{|\nabla^N x_3|} 
\end{align*}
for $s \in (0,\Lambda^{\frac{1}{4}}]$. On the other hand, applying Lemma \ref{technical.ingredient.2} with $\psi = e^{r_0 (1-H)} \, \frac{1}{|\nabla^N x_3|}$ gives 
\begin{align*} 
&\rho(s)^{\frac{10}{11}} \int_{N \cap \{x_3=s\}} e^{-\frac{\rho(s)^2 \, (x_1^2+x_2^2)}{4\tau}+\rho(s) \, (q_1x_1+q_2x_2)-r_0 H} \, \frac{1}{|\nabla^N x_3|} \\ 
&\leq \int_{N \cap \{x_3=s\}} e^{-\frac{x_1^2+x_2^2}{4\tau}+q_1x_1+q_2x_2-r_0 H} \, \frac{1}{|\nabla^N x_3|} 
\end{align*}
for $s \in (0,\Lambda^{\frac{1}{4}}]$. Putting these facts together, we conclude that 
\begin{align*} 
&\int_{\tilde{N} \cap \{x_3=s\}} e^{-\frac{x_1^2+x_2^2}{4\tau}+q_1x_1+q_2x_2-r_0 \tilde{H}} \, \frac{1}{|\nabla^{\tilde{N}} x_3|} \\ 
&\leq \rho(s)^{\frac{1}{11}} \int_{N \cap \{x_3=s\}} e^{-\frac{x_1^2+x_2^2}{4\tau}+q_1x_1+q_2x_2-r_0 H} \, \frac{1}{|\nabla^N x_3|} 
\end{align*} 
for $s \in (0,\Lambda^{\frac{1}{4}}]$. This proves the assertion. \\

We now consider the intermediate region $s \in (\Lambda^{\frac{1}{4}},\frac{\Lambda}{4}]$.

\begin{lemma}
\label{intermediate.region.a}
We can find positive real numbers $\delta_*$ and $\Lambda_*$, and a positive function $E_*(\Lambda)$ with the following property. If $\hat{\delta} < \delta_*$, $\Lambda > \Lambda_*$, and $\varepsilon < E_*(\Lambda)$, then we have 
\begin{align*} 
&\int_{\tilde{N} \cap \{x_3=s\}} e^{-\frac{x_1^2+x_2^2}{4\tau}+q_1x_1+q_2x_2-r_0 \tilde{H}} \, \frac{1}{|\nabla^{\tilde{N}} x_3|} \\ 
&\leq \int_{N \cap \{x_3=s\}} e^{-\frac{x_1^2+x_2^2}{4\tau}+q_1x_1+q_2x_2-r_0 H} \, \frac{1}{|\nabla^N x_3|} 
\end{align*}
for all $s \in (\Lambda^{\frac{1}{4}},\frac{\Lambda}{4}]$, all $r_0 \in [\frac{1}{1000},1]$, all $\tau \geq \frac{5}{9}$, and all $(q_1,q_2) \in \mathbb{R}^2$. 
\end{lemma}

\textbf{Proof.}
We first observe that the claim is true if $\sqrt{q_1^2+q_2^2}$ is sufficiently large. More precisely, we can find positive real numbers $\delta_*$ and $\Lambda_1$, and positive functions $E_1(\Lambda)$ and $Q(\Lambda)$ with the following significance. If $\hat{\delta} < \delta_*$, $\Lambda > \Lambda_1$, and $\varepsilon < E_1(\Lambda)$, then we have 
\begin{align*} 
&\int_{\tilde{N} \cap \{x_3=s\}} e^{-\frac{x_1^2+x_2^2}{4\tau}+q_1x_1+q_2x_2-r_0 \tilde{H}} \, \frac{1}{|\nabla^{\tilde{N}} x_3|} \\ 
&\leq \int_{N \cap \{x_3=s\}} e^{-\frac{x_1^2+x_2^2}{4\tau}+q_1x_1+q_2x_2-r_0 H} \, \frac{1}{|\nabla^N x_3|} 
\end{align*}
whenever $s \in (\Lambda^{\frac{1}{4}},\frac{\Lambda}{4}]$, $r_0 \in [\frac{1}{1000},1]$, $\tau \geq \frac{5}{9}$, and $\sqrt{q_1^2+q_2^2} > Q(\Lambda)$.

Therefore, it remains to consider the case $\sqrt{q_1^2+q_2^2} \leq Q(\Lambda)$. For each $s \in (\Lambda^{\frac{1}{4}},\frac{\Lambda}{4}]$, the cross section $\tilde{N} \cap \{x_3=s\}$ is obtained by dilating the cross section $N \cap \{x_3=s\}$ by the factor $\rho(s) := 1-e^{-\frac{4\Lambda}{s}}$, up to errors of order $O(\varepsilon)$. Moreover, we have the pointwise inequalities 
\[\frac{|\nabla^N x_3|}{|\nabla^{\tilde{N}} x_3|} 
\leq 1 + C \, \Big ( \frac{\Lambda}{x_3^2} \, e^{-\frac{4\Lambda}{x_3}} \Big )^2 + C \, \varepsilon \leq 1 + C \, \Lambda^{-2} \, e^{-\frac{4\Lambda}{x_3}} + C \, \varepsilon\] 
and 
\[\tilde{H} - H \geq 0\] 
for $x_3 \in (\Lambda^{\frac{1}{4}},\frac{\Lambda}{4}]$. Hence, we obtain 
\begin{align*} 
&\int_{\tilde{N} \cap \{x_3=s\}} e^{-\frac{x_1^2+x_2^2}{4\tau}+q_1x_1+q_2x_2-r_0 \tilde{H}} \, \frac{1}{|\nabla^{\tilde{N}} x_3|} \\ 
&\leq (1+C \, \Lambda^{-2} \, e^{-\frac{4\Lambda}{s}}+C(\Lambda) \, \varepsilon) \, \rho(s) \\ 
&\hspace{20mm} \cdot \int_{N \cap \{x_3=s\}} e^{-\frac{\rho(s)^2 \, (x_1^2+x_2^2)}{4\tau}+\rho(s) \, (q_1x_1+q_2x_2)-r_0 H} \, \frac{1}{|\nabla^N x_3|} 
\end{align*} 
whenever $s \in (\Lambda^{\frac{1}{4}},\frac{\Lambda}{4}]$, $r_0 \in [\frac{1}{1000},1]$, $\tau \geq \frac{5}{9}$, and $\sqrt{q_1^2+q_2^2} \leq Q(\Lambda)$. On the other hand, using Lemma \ref{technical.ingredient.2} with $\psi = e^{r_0 (1-H)} \, \frac{1}{|\nabla^N x_3|}$ gives 
\begin{align*} 
&\rho(s)^{\frac{10}{11}} \int_{N \cap \{x_3=s\}} e^{-\frac{\rho(s)^2 \, (x_1^2+x_2^2)}{4\tau}+\rho(s) \, (q_1x_1+q_2x_2)-r_0 H} \, \frac{1}{|\nabla^N x_3|} \\ 
&\leq \int_{N \cap \{x_3=s\}} e^{-\frac{x_1^2+x_2^2}{4\tau}+q_1x_1+q_2x_2-r_0 H} \, \frac{1}{|\nabla^N x_3|} 
\end{align*} 
whenever $s \in (\Lambda^{\frac{1}{4}},\frac{\Lambda}{4}]$, $r_0 \in [\frac{1}{1000},1]$, and $\tau \geq \frac{5}{9}$. Putting these facts together, we obtain 
\begin{align*} 
&\int_{\tilde{N} \cap \{x_3=s\}} e^{-\frac{x_1^2+x_2^2}{4\tau}+q_1x_1+q_2x_2-r_0 \tilde{H}} \, \frac{1}{|\nabla^{\tilde{N}} x_3|} \\ 
&\leq (1+C \, \Lambda^{-2} \, e^{-\frac{4\Lambda}{s}} + C(\Lambda) \, \varepsilon) \, \rho(s)^{\frac{1}{11}} \\ 
&\hspace{20mm} \cdot \int_{N \cap \{x_3=s\}} e^{-\frac{x_1^2+x_2^2}{4\tau}+q_1x_1+q_2x_2-r_0 H} \, \frac{1}{|\nabla^N x_3|} 
\end{align*} 
whenever $s \in (\Lambda^{\frac{1}{4}},\frac{\Lambda}{4}]$, $r_0 \in [\frac{1}{1000},1]$, $\tau \geq \frac{5}{9}$, and $\sqrt{q_1^2+q_2^2} \leq Q(\Lambda)$. Clearly, we can find a positive constant $\Lambda_2$ and a positive function $E_2(\Lambda)$ such that 
\[(1+C \, \Lambda^{-2} \, e^{-\frac{4\Lambda}{s}} + C(\Lambda) \, \varepsilon) \, \rho(s)^{\frac{1}{11}} \leq 1\] 
if $\Lambda > \Lambda_2$ and $\varepsilon < E_2(\Lambda)$. Hence, if we put $\Lambda_* = \max\{\Lambda_1,\Lambda_2\}$ and $E_*(\Lambda) = \min \{E_1(\Lambda),E_2(\Lambda)\}$, then the assertion follows. \\

In the next step, we consider the region $s \in (\frac{\Lambda}{4},\Lambda+\Lambda^{\frac{1}{4}}]$.

\begin{lemma}
\label{intermediate.region.b}
We can find positive real numbers $\delta_*$ and $\Lambda_*$ with the following property. If $\hat{\delta} < \delta_*$ and $\Lambda > \Lambda_*$, then we have 
\begin{align*} 
&\int_{\tilde{N} \cap \{x_3=s\}} e^{-\frac{x_1^2+x_2^2}{4\tau}+q_1x_1+q_2x_2-r_0 \tilde{H}} \, \frac{1}{|\nabla^{\tilde{N}} x_3|} \\ 
&\leq \int_{N \cap \{x_3=s\}} e^{-\frac{x_1^2+x_2^2}{4\tau}+q_1x_1+q_2x_2-r_0 H} \, \frac{1}{|\nabla^N x_3|} 
\end{align*}
for all $s \in (\frac{\Lambda}{4},\Lambda+\Lambda^{\frac{1}{4}}]$, all $r_0 \in [\frac{1}{1000},1]$, all $\tau \geq \frac{5}{9}$, and all $(q_1,q_2) \in \mathbb{R}^2$. 
\end{lemma}

\textbf{Proof.} 
Again, the claim is true if $\sqrt{q_1^2+q_2^2}$ is sufficiently large. More precisely, we can find positive real numbers $\delta_1$ and $\Lambda_1$, and a positive constant $Q$ with the following significance. If $\hat{\delta} < \delta_1$ and $\Lambda > \Lambda_1$, then we have 
\begin{align*} 
&\int_{\tilde{N} \cap \{x_3=s\}} e^{-\frac{x_1^2+x_2^2}{4\tau}+q_1x_1+q_2x_2-r_0 \tilde{H}} \, \frac{1}{|\nabla^{\tilde{N}} x_3|} \\ 
&\leq \int_{N \cap \{x_3=s\}} e^{-\frac{x_1^2+x_2^2}{4\tau}+q_1x_1+q_2x_2-r_0 H} \, \frac{1}{|\nabla^N x_3|} 
\end{align*}
whenever $s \in (\frac{\Lambda}{4},\Lambda+\Lambda^{\frac{1}{4}}]$, $r_0 \in [\frac{1}{1000},1]$, $\tau \geq \frac{5}{9}$, and $\sqrt{q_1^2+q_2^2} > Q$. 

Hence, it remains to consider the case $\sqrt{q_1^2+q_2^2} \leq Q$. For each $s \in (\frac{\Lambda}{4},\Lambda+\Lambda^{\frac{1}{4}}]$, the cross section $N \cap \{x_3=s\}$ is close to a circle of radius $1$, whereas $\tilde{N} \cap \{x_3=s\}$ is close to a circle of radius $\rho(s) := 1-e^{-\frac{4\Lambda}{s}}$. Moreover, $|\nabla^N x_3|$ and $|\nabla^{\tilde{N}} x_3|$ are close to $1$. By Lemma \ref{technical.ingredient.2}, we have 
\begin{align*} 
&\int_{\{x_1^2+x_2^2=\rho(s)^2\}} e^{-\frac{x_1^2+x_2^2}{4\tau}+q_1x_1+q_2x_2} \\ 
&\leq \rho(s)^{\frac{1}{11}} \int_{\{x_1^2+x_2^2=1\}} e^{-\frac{x_1^2+x_2^2}{4\tau}+q_1x_1+q_2x_2} 
\end{align*}
for all $(q_1,q_2) \in \mathbb{R}^2$. Hence, we can find positive real numbers $\delta_2$ and $\Lambda_2$ with the following property: if $\hat{\delta} < \delta_2$ and $\Lambda > \Lambda_2$, then 
\begin{align*} 
&\int_{\tilde{N} \cap \{x_3=s\}} e^{-\frac{x_1^2+x_2^2}{4\tau}+q_1x_1+q_2x_2-r_0 \tilde{H}} \, \frac{1}{|\nabla^{\tilde{N}} x_3|} \\ 
&\leq \int_{N \cap \{x_3=s\}} e^{-\frac{x_1^2+x_2^2}{4\tau}+q_1x_1+q_2x_2-r_0 H} \, \frac{1}{|\nabla^N x_3|} 
\end{align*}
whenever $s \in (\frac{\Lambda}{4},\Lambda+\Lambda^{\frac{1}{4}}]$, $r_0 \in [\frac{1}{1000},1]$, $\tau \geq \frac{5}{9}$, and $\sqrt{q_1^2+q_2^2} \leq Q$. Therefore, if we put $\delta_* = \min \{\delta_1,\delta_2\}$ and $\Lambda = \max \{\Lambda_1,\Lambda_2\}$, then the assertion follows. \\

Finally, we consider the region $s \in (\Lambda+\Lambda^{\frac{1}{4}},\Lambda+2 \, \Lambda^{\frac{1}{4}})$.

\begin{lemma}
\label{tip}
There exist positive real numbers $\delta_*$ and $\Lambda_*$ with the following property. If $\hat{\delta} < \delta_*$ and $\Lambda > \Lambda_*$, then we have 
\begin{align*} 
&\int_{\tilde{N} \cap \{x_3=s\}} e^{-\frac{x_1^2+x_2^2}{4\tau}+q_1x_1+q_2x_2-r_0 \tilde{H}} \, \frac{1}{|\nabla^{\tilde{N}} x_3|} \\ 
&\leq \int_{N \cap \{x_3=s\}} e^{-\frac{x_1^2+x_2^2}{4\tau}+q_1x_1+q_2x_2-r_0 H} \, \frac{1}{|\nabla^N x_3|} 
\end{align*}
for all $s \in (\Lambda+\Lambda^{\frac{1}{4}},\Lambda+2 \, \Lambda^{\frac{1}{4}})$, all $r_0 \in [\frac{1}{1000},1]$, all $\tau \geq \frac{5}{9}$, and all $(q_1,q_2) \in \mathbb{R}^2$. 
\end{lemma}

\textbf{Proof.}
For each $s \in (\Lambda+\Lambda^{\frac{1}{4}},\Lambda+2 \, \Lambda^{\frac{1}{4}})$, the cross section of $\tilde{N}$ at height $s$ is a circle of radius 
\[\rho(s) = a \, \sqrt{\frac{b-s}{a+b-s}},\] 
where $a = 1-e^{-4} + \frac{1}{3} \, (1-e^{-4})^2 \, \Lambda^{-\frac{1}{4}} < 1$ and $b = \Lambda + 2 \, \Lambda^{\frac{1}{4}}$. Note that  
\[\frac{1}{|\nabla^{\tilde{N}} x_3|} = \sqrt{1+\rho'(s)^2}\] 
and 
\[\tilde{H} \geq \frac{1}{a}\] 
on the set $\tilde{N} \cap \{x_3=s\}$. Therefore, 
\begin{align*} 
&\int_{\tilde{N} \cap \{x_3=s\}} e^{-\frac{x_1^2+x_2^2}{4\tau}+q_1x_1+q_2x_2-r_0 \tilde{H}} \, \frac{1}{|\nabla^{\tilde{N}} x_3|} \\ 
&\leq \int_{\tilde{N} \cap \{x_3=s\}} e^{-\frac{x_1^2+x_2^2}{4\tau}+q_1x_1+q_2x_2-\frac{r_0}{a}} \, \frac{1}{|\nabla^{\tilde{N}} x_3|} \\ 
&= \rho(s) \, \sqrt{1+\rho'(s)^2} \, e^{-\frac{\rho(s)^2}{4\tau}-\frac{r_0}{a}} \int_{S^1} e^{\rho(s) \, (q_1x_1+q_2x_2)} \\ 
&= a \, \sqrt{1-\frac{a}{a+b-s}+\frac{a^4}{4 \, (a+b-s)^4}} \, e^{\frac{a^2}{4\tau} \, \frac{a}{a+b-s}-\frac{a^2}{4\tau}-\frac{r_0}{a}} \\ 
&\hspace{10mm} \cdot \int_{S^1} e^{a \, \sqrt{\frac{b-s}{a+b-s}} \, (q_1x_1+q_2x_2)}
\end{align*} 
for all $\tau \geq \frac{5}{9}$ and all $s \in (\Lambda+\Lambda^{\frac{1}{4}},\Lambda+2 \, \Lambda^{\frac{1}{4}})$. Here, $S^1$ denotes the unit circle of radius $1$. It is elementary to check that 
\[1-z+\frac{z^4}{4} \leq e^{-z},\] 
hence 
\[\sqrt{1-z+\frac{z^4}{4}} \, e^{\frac{z}{2}} \leq 1\] 
for all $z \in [0,1]$. Consequently, we have 
\[\sqrt{1-\frac{a}{a+b-s}+\frac{a^4}{4 \, (a+b-s)^4}} \, e^{\frac{a^2}{4\tau} \, \frac{a}{a+b-s}} \leq 1\] 
for all $\tau \geq \frac{5}{9}$ and $s \in (\Lambda+\Lambda^{\frac{1}{4}},\Lambda+2 \, \Lambda^{\frac{1}{4}})$. From this, we deduce that 
\begin{align*} 
&\int_{\tilde{N} \cap \{x_3=s\}} e^{-\frac{x_1^2+x_2^2}{4\tau}+q_1x_1+q_2x_2-r_0 \tilde{H}} \, \frac{1}{|\nabla^{\tilde{N}} x_3|} \\ 
&\leq a \, e^{-\frac{a^2}{4\tau}-\frac{r_0}{a}} \, \int_{S^1} e^{a \, \sqrt{\frac{b-s}{a+b-s}} \, (q_1x_1+q_2x_2)} 
\end{align*} 
for all $\tau \geq \frac{5}{9}$ and $s \in (\Lambda+\Lambda^{\frac{1}{4}},\Lambda+2 \, \Lambda^{\frac{1}{4}})$. On the other hand, if the neck $N$ is sufficiently close to a cylinder of radius $1$, then we have 
\begin{align*} 
&\int_{S^1} e^{a \, \sqrt{\frac{b-s}{a+b-s}} \, (q_1x_1+q_2x_2)} \\ 
&\leq \int_{N \cap \{x_3=s\}} e^{\frac{1-x_1^2-x_2^2}{4\tau}+q_1x_1+q_2x_2+r_0 (\frac{1}{a}-H)} \, \frac{1}{|\nabla^N x_3|} 
\end{align*} 
for all $s \in (\Lambda+\Lambda^{\frac{1}{4}},\Lambda+2 \, \Lambda^{\frac{1}{4}})$, all $r_0 \in [\frac{1}{1000},1]$, all $\tau \geq \frac{5}{9}$, and all $(q_1,q_2) \in \mathbb{R}^2$. Putting these facts together, we conclude that 
\begin{align*} 
&\int_{\tilde{N} \cap \{x_3=s\}} e^{-\frac{x_1^2+x_2^2}{4\tau}+q_1x_1+q_2x_2-r_0 \tilde{H}} \, \frac{1}{|\nabla^{\tilde{N}} x_3|} \\ 
&\leq a \, e^{\frac{1-a^2}{4\tau}} \, 
\int_{N \cap \{x_3=s\}} e^{-\frac{x_1^2+x_2^2}{4\tau}+q_1x_1+q_2x_2-r_0 H} \, \frac{1}{|\nabla^N x_3|} 
\end{align*} 
for all $s \in (\Lambda+\Lambda^{\frac{1}{4}},\Lambda+2 \, \Lambda^{\frac{1}{4}})$, all $r_0 \in [\frac{1}{1000},1]$, all $\tau \geq \frac{5}{9}$, and all $(q_1,q_2) \in \mathbb{R}^2$. Since $a$ is close to $1-e^{-4}$ for $\Lambda$ large, we have $a \, e^{\frac{1-a^2}{4\tau}} \leq a \, e^{\frac{9(1-a^2)}{20}} \leq 1$. From this, the assertion follows. \\

Combining the previous results, we can draw the following conclusion: 

\begin{proposition}
\label{monotonicity.during.single.surgery}
There exist real numbers $\delta_*$ and $\Lambda_*$, and a function $E_*(\Lambda)$ with the following property. If $\hat{\delta} < \delta_*$, $\Lambda > \Lambda_*$, and $\varepsilon < E_*(\Lambda)$, then we have 
\[\int_{\tilde{N} \cap \{0 \leq x_3 \leq \Lambda+2\Lambda^{\frac{1}{4}}\}} e^{-\frac{|x-p|^2}{4\tau}-r_0 \tilde{H}} \leq \int_{N \cap \{0 \leq x_3 \leq \Lambda+2\Lambda^{\frac{1}{4}}\}} e^{-\frac{|x-p|^2}{4\tau}-r_0 H}\] 
for all $r_0 \in [\frac{1}{1000},1]$, all $\tau \geq \frac{5}{9}$, and all $p \in \mathbb{R}^3$.
\end{proposition}

\textbf{Proof.} 
In view of Lemma \ref{x_3.small}, Lemma \ref{intermediate.region.a}, Lemma \ref{intermediate.region.b}, and Lemma \ref{tip}, we can find positive real numbers $\delta_*$, and $\Lambda_*$, and a positive function $E_*(\Lambda)$ with the following property. If $\hat{\delta} < \delta_*$, $\Lambda > \Lambda_*$, and $\varepsilon < E_*(\Lambda)$, then we have 
\begin{align*} 
&\int_{\tilde{N} \cap \{x_3=s\}} e^{-\frac{x_1^2+x_2^2}{4\tau}+q_1x_1+q_2x_2-r_0 \tilde{H}} \, \frac{1}{|\nabla^{\tilde{N}} x_3|} \\ 
&\leq \int_{N \cap \{x_3=s\}} e^{-\frac{x_1^2+x_2^2}{4\tau}+q_1x_1+q_2x_2-r_0 H} \, \frac{1}{|\nabla^N x_3|} 
\end{align*}
for all $s \in (0,\Lambda+2\Lambda^{\frac{1}{4}})$, all $r_0 \in [\frac{1}{1000},1]$, all $\tau \geq \frac{5}{9}$, and all $(q_1,q_2) \in \mathbb{R}^2$. This implies 
\begin{align*} 
&\int_{\tilde{N} \cap \{x_3=s\}} e^{-\frac{|x-p|^2}{4\tau}-r_0 \tilde{H}} \, \frac{1}{|\nabla^{\tilde{N}} x_3|} \\ 
&\leq \int_{N \cap \{x_3=s\}} e^{-\frac{|x-p|^2}{4\tau}-r_0 H} \, \frac{1}{|\nabla^N x_3|} 
\end{align*}
for all $s \in (0,\Lambda+2\Lambda^{\frac{1}{4}})$, all $r_0 \in [\frac{1}{1000},1]$, all $\tau \geq \frac{5}{9}$, and all $p \in \mathbb{R}^3$. If we integrate over $s$ and apply the co-area formula, the assertion follows.

\section{A monotonicity formula for mean curvature flow with surgery in $\mathbb{R}^3$}

\label{monotonicity.for.mcf.with.surgery}

\begin{proposition}
\label{smooth.case}
Let $M_t$ be a family of mean convex surfaces in $\mathbb{R}^3$ which evolve under smooth mean curvature flow. For each $r_0>0$, the function 
\[\int_{M_t} \frac{1}{4\pi (t_0-t)} \, e^{-\frac{|x-p|^2}{4(t_0-t)}-r_0 H}\] 
is monotone decreasing for $t_0-t>0$.
\end{proposition}

\textbf{Proof.} 
We compute 
\begin{align*} 
\Big ( \frac{\partial}{\partial t}-\Delta \Big ) e^{-r_0 H} 
&= -r_0 \, e^{-r_0 H} \, \Big ( \frac{\partial}{\partial t}-\Delta \Big ) H - r_0^2 \, e^{-r_0 H} \, |\nabla H|^2 \\ 
&= -r_0 \, e^{-r_0 H} \, |A|^2 \, H - r_0^2 \, e^{-r_0 H} \, |\nabla H|^2 \leq 0. 
\end{align*} 
Hence, the assertion follows from Ecker's weighted monotonicity formula (see Theorem 4.13 in \cite{Ecker}). \\

We next consider a mean curvature flow with surgery in $\mathbb{R}^3$. We assume that each surgery procedure involves performing $\Lambda$-surgery on an $(\hat{\alpha},\hat{\delta},\varepsilon,L)$-neck of size $r \in [\frac{1}{2H_1},\frac{1}{H_1}]$ (see \cite{Brendle-Huisken1} for definitions).

\begin{theorem}
\label{monotonicity.euclidean.case}
Let $\delta_*$, $\Lambda_*$, and $E_*(\Lambda)$ be defined as in Proposition \ref{monotonicity.during.single.surgery}. Moreover, suppose that $M_t$ is a mean curvature flow with surgery, and that the surgery parameters satisfy $\hat{\delta} < \delta_*$, $\Lambda > \Lambda_*$, and $\varepsilon < E_*(\Lambda)$. Then the function 
\[\int_{M_t} \frac{1}{4\pi (t_0-t)} \, e^{-\frac{|x-p|^2}{4(t_0-t)}-\frac{H}{200 \, H_1}}\] 
is monotone decreasing for $t_0-t \geq \frac{5}{9} \, H_1^{-2}$. 
\end{theorem}

\textbf{Proof.} 
Proposition \ref{smooth.case} guarantees that the monotonicity formula holds in between surgery times. Moreover, it follows from Proposition \ref{monotonicity.during.single.surgery} that the monotonicity property holds across surgery times. \\

Theorem \ref{monotonicity.euclidean.case} allows us to draw the following conclusion: 

\begin{corollary}
\label{application}
We can find positive real numbers $\delta_*$, $\Lambda_*$, $L$, and a positive function $E_*(\Lambda)$ with the following property. Suppose that $M_t$ is a mean curvature flow with surgery, and that the surgery parameters satisfy $\hat{\delta} < \delta_*$, $\Lambda > \Lambda_*$, and $\varepsilon < E_*(\Lambda)$. Finally, suppose that $M_{t_0-}$ contains an $(\hat{\alpha},\hat{\delta},\varepsilon,L)$-neck of size $r \in [\frac{1}{2H_1},\frac{1}{H_1}]$, and $p$ is a point in ambient space which lies on the axis of that neck. Then 
\[\int_{M_{t_0+\frac{5}{9} \, H_1^{-2}-\tau}} \frac{1}{4\pi \tau} \, e^{-\frac{|x-p|^2}{4\tau}} \geq 1.01\] 
for all $\tau \geq \frac{5}{9} \, H_1^{-2}$.
\end{corollary} 

\textbf{Proof.} 
It is elementary to check that 
\[\inf_{r \in [\frac{1}{2H_1},\frac{1}{H_1}]} H_1 \, r \, e^{-\frac{9 \, H_1^2 \, r^2}{20}-\frac{1}{200 \, H_1 \, r}} = \frac{1}{2} \, e^{-\frac{49}{400}}.\] 
Hence, on an exact cylinder $S^1(r) \times \mathbb{R}$, we have 
\begin{align*} 
\int_{S^1(r) \times \mathbb{R}} \frac{9 \, H_1^2}{20\pi} \, e^{-\frac{9 \, H_1^2 \, |x|^2}{20}-\frac{H}{200 \, H_1}} 
&= \sqrt{\frac{9\pi}{5}} \, H_1 \, r \, e^{-\frac{9 \, H_1^2 \, r^2}{20}-\frac{1}{200 \, H_1 \, r}} \\ 
&\geq \sqrt{\frac{9\pi}{20}} \, e^{-\frac{49}{400}} \\ 
&\geq 1.02 
\end{align*}
for all $r \in [\frac{1}{2H_1},\frac{1}{H_1}]$. By assumption, the surface $M_{t_0-}$ contains a neck of size $r \in [\frac{1}{2H_1},\frac{1}{H_1}]$, and $p$ lies on the axis of that neck. Hence, if $L$ is sufficiently large, then we obtain 
\[\int_{M_{t_0-}} \frac{9 \, H_1^2}{20\pi} \, e^{-\frac{9 \, H_1^2 \, |x-p|^2}{20}-\frac{H}{200 \, H_1}} \geq 1.01.\] 
Using Theorem \ref{monotonicity.euclidean.case}, we obtain 
\begin{align*} 
&\int_{M_{t_0+\frac{5}{9} \, H_1^{-2}-\tau}} \frac{1}{4\pi \tau} \, e^{-\frac{|x-p|^2}{4\tau}} \\ 
&\geq \int_{M_{t_0+\frac{5}{9} \, H_1^{-2}-\tau}} \frac{1}{4\pi \tau} \, e^{-\frac{|x-p|^2}{4\tau}-\frac{H}{200 \, H_1}} \\ 
&\geq \int_{M_{t_0-}} \frac{9 \, H_1^2}{20\pi} \, e^{-\frac{9 \, H_1^2 \, |x-p|^2}{20}-\frac{H}{200 \, H_1}} \\ 
&\geq 1.01
\end{align*} 
provided that $\tau \geq \frac{5}{9} \, H_1^{-2}$ and $L$ is sufficiently large. \\

\begin{theorem}
\label{regularity}
Fix an open interval $I$ and a compact interval $J \subset I$. Suppose that $\bar{\mathcal{M}}$ is an embedded smooth solution of mean curvature flow which is defined for $t \in I$. Moreover, suppose that $\mathcal{M}^{(j)}$ is a sequence of mean curvature flows with surgery, each of which is defined for $t \in I$. We assume that each surgery of the flow $\mathcal{M}^{(j)}$ involves performing $\Lambda$-surgery on an $(\hat{\alpha},\hat{\delta},\varepsilon,L)$-neck of size $r \in [\frac{1}{2H_1^{(j)}},\frac{1}{H_1^{(j)}}]$, where $H_1^{(j)} \to \infty$. We assume further that the surgery parameters satisfy $\hat{\delta} < \delta_*$, $\Lambda > \Lambda_*$, and $\varepsilon < E_*(\Lambda)$. Finally, we assume that the flows $\mathcal{M}^{(j)}$ converge to $\bar{\mathcal{M}}$ in the sense of geometric measure theory. Then, if $j$ is sufficiently large, the flow $\mathcal{M}^{(j)}$ is free of surgeries for all $t \in J$. Furthermore, the flows $\mathcal{M}^{(j)}$ converge smoothly to $\bar{\mathcal{M}}$ as $j \to \infty$. 
\end{theorem}

\textbf{Proof.}
We first show that, for $j$ sufficiently large, the flow $\mathcal{M}^{(j)}$ is free of surgeries for all $t \in J$. Suppose that each flow $\mathcal{M}^{(j)}$ has at least one surgery time $t_j \in J$. For each $j$, we can find an $(\hat{\alpha},\hat{\delta},\varepsilon,L)$-neck in $M_{t_j-}^{(j)}$ of size $r_j \in [\frac{1}{2H_1^{(j)}},\frac{1}{H_1^{(j)}}]$. Let $p_j$ be a point in ambient space which lies on the axis of that neck. Using Corollary \ref{application}, we obtain 
\[\int_{M_{t_j+\frac{5}{9} \, (H_1^{(j)})^{-2}-\tau}^{(j)}} \frac{1}{4\pi \tau} \, e^{-\frac{|x-p_j|^2}{4\tau}} \geq 1.01\] 
for all $\tau \geq \frac{5}{9} \, (H_1^{(j)})^{-2}$. We now pass to the limit as $j \to \infty$. If we define $\bar{t} = \lim_{j \to \infty} t_j \in J$ and $\bar{p} = \lim_{j \to \infty} p_j$, then we obtain 
\[\int_{\bar{M}_{\bar{t}-\tau}} \frac{1}{4\pi \tau} \, e^{-\frac{|x-\bar{p}|^2}{4\tau}} \geq 1.01\] 
for each $\tau \in (0,\inf J - \inf I)$. On the other hand, since $\bar{\mathcal{M}}$ is smooth, we have 
\[\int_{\bar{M}_{\bar{t}-\tau}} \frac{1}{4\pi \tau} \, e^{-\frac{|x-\bar{p}|^2}{4\tau}} \to 1\] 
as $\tau \to 0$. This is a contradiction. Therefore, the flow $\mathcal{M}^{(j)}$ is free of surgeries for $t \in J$. Using standard local regularity theorems for mean curvature flow (cf. \cite{Brakke},\cite{White3}), we conclude that the flows $\mathcal{M}^{(j)}$ converge smoothly to $\bar{\mathcal{M}}$ as $j \to \infty$. 

\section{Adaptation to the Riemannian setting}

\label{Riem.manifolds}

Let $M_t$ be a mean convex solution of mean curvature flow in a compact Riemannian three-manifold $X$. Let us fix a time $t_0$ and point $p$ in ambient space. Let $\varphi$ be a smooth cutoff on $X$ such that $\varphi(x) = 1$ for $d(p,x) \leq \frac{1}{4} \, \text{\rm inj}(X)$ and $\varphi(x) = 0$ for $d(p,x) \geq \frac{1}{2} \, \text{\rm inj}(X)$. Moreover, we put 
\[\Phi(x,t) = \frac{1}{4\pi (t_0-t)} \, e^{-\frac{d(p,x)^2}{4(t_0-t)}} \, \varphi(x)^2\] 
for $t_0-t > 0$. 

\begin{proposition} 
\label{monotonicity.formula.in.manifolds} 
Let $K$ be a positive constant with the property that the ambient three-manifold $X$ has Ricci curvature is at least $-K$. Then we have 
\[\frac{d}{dt} \bigg ( \int_{M_t} \Phi \, \exp \big ( (t_0-t)^{\frac{1}{2}} - r_0 \, e^{-K \, (t_0-t)} \, H \big ) \bigg ) \leq C \, |M_t|\] 
whenever $t_0-t \in (0,1]$. The constant $C$ depends only on the ambient three-manifold $X$.
\end{proposition}

\textbf{Proof.} 
In the following, we assume that $t_0-t \in (0,1]$. In the region $d(p,x) \leq \frac{1}{4} \, \text{\rm inj}(X)$, we have 
\[\Big ( \frac{\partial}{\partial t} + \Delta - H^2 \Big ) \Phi 
\leq -\Big ( H + \frac{\bar{D}_\nu \Phi}{\Phi} \Big )^2 \, \Phi + C \, \frac{d(p,x)^2}{t_0-t} \, \Phi.\] 
This implies 
\[\Big ( \frac{\partial}{\partial t} + \Delta - H^2 \Big ) \Phi 
\leq C \, \frac{d(p,x)^2}{t_0-t} \, \Phi,\] 
hence 
\begin{align*} 
&\Big ( \frac{\partial}{\partial t} + \Delta - H^2 \Big ) \big ( \Phi \, \exp \big ( (t_0-t)^{\frac{1}{2}} \big ) \big ) \\ 
&\leq \Big ( C \, \frac{d(p,x)^2}{t_0-t} - \frac{1}{2 \, (t_0-t)^{\frac{1}{2}}} \Big ) \, \Phi \, \exp \big ( (t_0-t)^{\frac{1}{2}} \big ) \leq C 
\end{align*}
for $d(p,x) \leq \frac{1}{4} \, \text{\rm inj}(X)$. Moreover, in the region $\frac{1}{4} \, \text{\rm inj}(X) \leq d(p,x) \leq \frac{1}{2} \, \text{\rm inj}(X)$, we have 
\[\Big ( \frac{\partial}{\partial t} + \Delta - H^2 \Big ) \Phi \leq -\Big ( H + \frac{\bar{D}_\nu \Phi}{\Phi} \Big )^2 \, \Phi + C.\] 
This gives 
\[\Big ( \frac{\partial}{\partial t} + \Delta - H^2 \Big ) \big ( \Phi \, \exp \big ( (t_0-t)^{\frac{1}{2}} \big ) \big ) \leq C\] 
for $\frac{1}{4} \, \text{\rm inj}(X) \leq d(p,x) \leq \frac{1}{2} \, \text{\rm inj}(X)$. To summarize, we have shown that 
\[\Big ( \frac{\partial}{\partial t} + \Delta - H^2 \Big ) \big ( \Phi \, \exp \big ( (t_0-t)^{\frac{1}{2}} \big ) \big ) \leq C\] 
at each point in $M_t$. Since the Ricci curvature of $X$ is bounded from below by $-K$, we have 
\[\Big ( \frac{\partial}{\partial t} - \Delta \Big ) (e^{-K (t_0-t)} \, H) \geq e^{-K (t_0-t)} \, |A|^2 \, H \geq 0,\] 
hence 
\begin{align*} 
&\Big ( \frac{\partial}{\partial t}-\Delta \Big ) \exp \big ( -r_0 \, e^{-K (t_0-t)} \, H \big ) \\ 
&= -\exp \big ( -r_0 \, e^{-K (t_0-t)} \, H \big ) \, \Big ( \frac{\partial}{\partial t}-\Delta \Big ) (r_0 \, e^{-K (t_0-t)} \, H) \\ 
&- \exp \big ( -r_0 \, e^{-K (t_0-t)} \, H \big ) \, |\nabla (r_0 \, e^{-K (t_0-t)} \, H)|^2 \\ 
&\leq 0. 
\end{align*}
Putting these facts together, we obtain 
\begin{align*} 
&\frac{d}{dt} \bigg ( \int_{M_t} \Phi \, \exp \big ( (t_0-t)^{\frac{1}{2}} - r_0 \, e^{-K (t_0-t)} \, H \big ) \bigg ) \\ 
&\leq C \int_{M_t} \exp \big ( -r_0 \, e^{-K (t_0-t)} \, H \big ) \leq C \, |M_t|, 
\end{align*}
provided that $t_0-t \in (0,1]$. This completes the proof. \\

We now consider mean curvature flow with surgery in a Riemannian manifold. We begin with a definition.

\begin{definition}
Let $M$ be a closed surface in a Riemannian three-manifold, and let $N$ be a region in $M$. We say that $N$ is an $(\hat{\alpha},\hat{\delta},\varepsilon,L)$-neck of size $r$ if there exists a point $o \in N$ with the property that the surface $\exp_o^{-1}(N)$ is an $(\hat{\alpha},\hat{\delta},\varepsilon,L)$-neck of size $r$ in Euclidean space (see \cite{Brendle-Huisken1} for the definition).
\end{definition}

We next explain how to do on a neck in Riemannian three-manifold. Namely, if $o$ lies at the center of an $(\hat{\alpha},\hat{\delta},\varepsilon,L)$-neck $N$ in $X$, then $\exp_o^{-1}(N)$ is an $(\hat{\alpha},\hat{\delta},\varepsilon,L)$-neck in Euclidean space. Hence, we can perform the surgery procedure described in \cite{Brendle-Huisken1} on $\exp_o^{-1}(N)$. We then paste the surgically modified surface back into $X$ using the exponential map $\exp_o$.

In the next step, we analyze how the quantity in Proposition \ref{monotonicity.formula.in.manifolds} changes under a single surgery. To that end, it will be convenient to work in geodesic normal coordinates around $o$; that is, we will identify a point in $T_o N = \mathbb{R}^3$ with its image under the exponential map $\exp_o$. With this identification, we can view $N$ as an $(\hat{\alpha},\hat{\delta},\varepsilon,L)$-neck in $\mathbb{R}^3$. Without loss of generality, we may assume that the axis of the neck $N$ is parallel to the $x_3$-axis. Note that the origin lies on $N$, so the axis of the neck does not pass through the origin. Finally, we denote by $\tilde{N}$ the surface obtained from $N$ by performing a $\Lambda$-surgery on $N$.

\begin{proposition} 
\label{monotonicity.during.single.surgery.riemannian.case}
Suppose that the surgery parameters $\hat{\delta}$, $\Lambda$, and $\varepsilon$ satisfy $\hat{\delta} < \delta_*$, $\Lambda > \Lambda_*$, and $\varepsilon < E_*(\Lambda)$, where $\delta_*$, $\Lambda_*$, and $E_*$ are defined as in Proposition \ref{monotonicity.during.single.surgery}. Moreover, let $N$ be an $(\hat{\alpha},\hat{\delta},\varepsilon,L)$-neck $N$ of size $r$ in $X$, and let $\tilde{N}$ denote the surgically modified surface. Then 
\begin{align*} 
&\int_{\tilde{N} \cap \{0 \leq x_3 \leq \Lambda+2\Lambda^{\frac{1}{4}}\}} e^{-\frac{d(p,x)^2}{4\tau} - r_0 \, e^{-K \tau} \, \tilde{H}} \, d\mu \\ 
&\leq \int_{N \cap \{0 \leq x_3 \leq \Lambda+2\Lambda^{\frac{1}{4}}\}} e^{-\frac{d(p,x)^2}{4\tau} - r_0 \, e^{-K \tau} \, H} \, d\mu + C(\Lambda) \, \tau \, r^2 
\end{align*}
provided that $d(o,p) \leq \frac{1}{8} \, \text{\rm inj}(X)$, $\frac{1}{K} \geq \tau \geq \frac{5}{9} \, r^2$ and $r_0 \in [\frac{1}{300} \, r,r]$.
\end{proposition}

\textbf{Proof.} 
It follows from Proposition \ref{monotonicity.during.single.surgery} that 
\begin{align*} 
&\int_{\tilde{N} \cap \{0 \leq x_3 \leq \Lambda+2\Lambda^{\frac{1}{4}}\}} e^{-\frac{|x-p|^2}{4\tau}-r_0 \, e^{-K \tau} \, \tilde{H}_{\text{\rm eucl}}} \, d\mu_{\text{\rm eucl}} \\ 
&\leq \int_{N \cap \{0 \leq x_3 \leq \Lambda+2\Lambda^{\frac{1}{4}}\}} e^{-\frac{|x-p|^2}{4\tau}-r_0 \, e^{-K \tau} \, H_{\text{\rm eucl}}} \, d\mu_{\text{\rm eucl}} 
\end{align*}
provided that $\tau \geq \frac{5}{9} \, r^2$ and $r_0 \in [\frac{1}{300} \, r,r]$. Here, $|x-p|$ denotes the Euclidean distance of $x$ and $p$ in geodesic normal coordinates around $o$; $H_{\text{\rm eucl}}$ and $\tilde{H}_{\text{\rm eucl}}$ denote the mean curvatures of $N$ and $\tilde{N}$ with respect to the Euclidean metric; and $d\mu_{\text{\rm eucl}}$ denotes the area form with respect to the Euclidean metric on $\mathbb{R}^3$. 

In the next step, we compare the Riemannian distance $d(p,x)$ to the Euclidean distance $|x-p|$. To that end, we perform a Taylor expansion of the function $x \mapsto \frac{1}{2} \, d(p,x)^2$ around the origin $o$. The value of this function at $o$ is given by $\frac{1}{2} \, |p|^2$. Its gradient at $o$ equals $-p$. Moreover, the Hessian of the function $x \mapsto \frac{1}{2} \, d(p,x)^2$ at $o$ equals $g + O(|p|^2)$. Finally, the third derivatives of the function $x \mapsto \frac{1}{2} \, d(p,x)^2$ at $o$ are bounded by $O(|p|)$. Putting these facts together, we obtain 
\[\frac{1}{2} \, d(p,x)^2 = \frac{1}{2} \, |p|^2 - \langle p,x \rangle + \frac{1}{2} \, |x|^2 + O(|p|^2 \, |x|^2 + |p| \, |x|^3 + |x|^4).\] 
In other words, 
\[\frac{1}{2} \, d(p,x)^2 = \frac{1}{2} \, |x-p|^2 + O((|p|^2+|x|^2) \, |x|^2).\] 
If $x \in (N \cup \tilde{N}) \cap \{0 \leq x_3 \leq \Lambda+2\Lambda^{\frac{1}{4}}\}$, then we have $|x|^2 \leq C(\Lambda) \, r^2$, hence $|x|^2 \leq C(\Lambda) \, \tau$. This implies 
\[\Big | \frac{d(p,x)^2}{4\tau} - \frac{|x-p|^2}{4\tau} \Big | \leq C(\Lambda) \, (|p|^2 + \tau)\] 
for all points $x \in (N \cup \tilde{N}) \cap \{0 \leq x_3 \leq \Lambda+2\Lambda^{\frac{1}{4}}\}$. Consequently, we have 
\[e^{-\frac{|x-p|^2}{4\tau}} \leq (1 + C(\Lambda) \, (|p|^2 + \tau)) \, e^{-\frac{d(p,x)^2}{4\tau}}\] 
for all points $x \in N \cap \{0 \leq x_3 \leq \Lambda+2\Lambda^{\frac{1}{4}}\}$. Since $|H-H_{\text{\rm eucl}}| \leq C(\Lambda) \, r$ for all points $x \in N \cap \{0 \leq x_3 \leq \Lambda+2\Lambda^{\frac{1}{4}}\}$, it follows that 
\[e^{-\frac{|x-p|^2}{4\tau} - r_0 \, e^{-K \tau} \, H_{\text{\rm eucl}}} \leq (1 + C(\Lambda) \, (|p|^2 + \tau)) \, e^{-\frac{d(p,x)^2}{4\tau} - r_0 \, e^{-K \tau} \, H}\] 
for all points $x \in N \cap \{0 \leq x_3 \leq \Lambda+2\Lambda^{\frac{1}{4}}\}$. Thus, 
\begin{align*} 
&\int_{N \cap \{0 \leq x_3 \leq \Lambda+2\Lambda^{\frac{1}{4}}\}} e^{-\frac{|x-p|^2}{4\tau} - r_0 \, e^{-K \tau} \, H_{\text{\rm eucl}}} \, d\mu_{\text{\rm eucl}} \\ 
&\leq (1 + C(\Lambda) \, (|p|^2 + \tau)) \int_{N \cap \{0 \leq x_3 \leq \Lambda+2\Lambda^{\frac{1}{4}}\}} e^{-\frac{d(p,x)^2}{4\tau} - r_0 \, e^{-K \tau} \, H} \, d\mu. 
\end{align*} 
Similarly, we have 
\[e^{-\frac{d(p,x)^2}{4\tau}} \leq (1 + C(\Lambda) \, (|p|^2 + \tau)) \, e^{-\frac{|x-p|^2}{4\tau}}\] 
for all points $x \in \tilde{N} \cap \{0 \leq x_3 \leq \Lambda+2\Lambda^{\frac{1}{4}}\}$. Since $|\tilde{H}-\tilde{H}_{\text{\rm eucl}}| \leq C(\Lambda) \, r$ for all points $x \in \tilde{N} \cap \{0 \leq x_3 \leq \Lambda+2\Lambda^{\frac{1}{4}}\}$, it follows that 
\[e^{-\frac{d(p,x)^2}{4\tau} - r_0 \, e^{-K \tau} \, \tilde{H}} \leq (1 + C(\Lambda) \, (|p|^2 + \tau)) \, e^{-\frac{|x-p|^2}{4\tau} - r_0 \, e^{-K \tau} \, \tilde{H}_{\text{\rm eucl}}}\] 
for all points $x \in \tilde{N} \cap \{0 \leq x_3 \leq \Lambda+2\Lambda^{\frac{1}{4}}\}$. This gives 
\begin{align*} 
&\int_{\tilde{N} \cap \{0 \leq x_3 \leq \Lambda+2\Lambda^{\frac{1}{4}}\}} e^{-\frac{d(p,x)^2}{4\tau} - r_0 \, e^{-K \tau} \, \tilde{H}} \, d\mu \\ 
&\leq (1 + C(\Lambda) \, (|p|^2 + \tau)) \int_{\tilde{N} \cap \{0 \leq x_3 \leq \Lambda+2\Lambda^{\frac{1}{4}}\}} e^{-\frac{|x-p|^2}{4\tau} - r_0 \, e^{-K \tau} \, \tilde{H}_{\text{\rm eucl}}} \, d\mu_{\text{\rm eucl}}. 
\end{align*} 
Putting these facts together, we conclude that 
\begin{align*} 
&\int_{\tilde{N} \cap \{0 \leq x_3 \leq \Lambda+2\Lambda^{\frac{1}{4}}\}} e^{-\frac{d(p,x)^2}{4\tau} - r_0 \, e^{-K \tau} \, \tilde{H}} \, d\mu \\ 
&\leq (1+C(\Lambda) \, (|p|^2+\tau)) \int_{N \cap \{0 \leq x_3 \leq \Lambda+2\Lambda^{\frac{1}{4}}\}} e^{-\frac{d(p,x)^2}{4\tau} - r_0 \, e^{-K \tau} \, H} \, d\mu. 
\end{align*}
Finally, we have the pointwise estimate $(|p|^2+\tau) \, e^{-\frac{d(p,x)^2}{4\tau}} \leq C(\Lambda) \, \tau$ at each point on $N \cap \{0 \leq x_3 \leq \Lambda+2\Lambda^{\frac{1}{4}}\}$. Therefore, we obtain 
\[(|p|^2+\tau) \int_{N \cap \{0 \leq x_3 \leq \Lambda+2\Lambda^{\frac{1}{4}}\}} e^{-\frac{d(p,x)^2}{4\tau} - r_0 \, e^{-K \tau} \, H} \, d\mu \leq C(\Lambda) \, \tau \, r^2,\] 
hence 
\begin{align*} 
&\int_{\tilde{N} \cap \{0 \leq x_3 \leq \Lambda+2\Lambda^{\frac{1}{4}}\}} e^{-\frac{d(p,x)^2}{4\tau} - r_0 \, e^{-K \tau} \, \tilde{H}} \, d\mu \\ 
&\leq \int_{N \cap \{0 \leq x_3 \leq \Lambda+2\Lambda^{\frac{1}{4}}\}} e^{-\frac{d(p,x)^2}{4\tau} - r_0 \, e^{-K \tau} \, H} \, d\mu + C(\Lambda) \, \tau \, r^2. 
\end{align*}
This proves the assertion. \\

Combining Proposition \ref{monotonicity.formula.in.manifolds} and Proposition \ref{monotonicity.during.single.surgery.riemannian.case}, we arrive at the following conclusion:

\begin{theorem}
\label{monotonicity.riemannian.case}
Let $M_t$ be a mean curvature flow with surgery in a Riemannian three-manifold $X$. Suppose that the surgery parameters satisfy $\hat{\delta} < \delta_*$, $\Lambda > \Lambda_*$, and $\varepsilon < E_*(\Lambda)$, where $\delta_*$, $\Lambda_*$, and $E_*(\Lambda)$ are defined as in Proposition \ref{monotonicity.during.single.surgery}. Then 
\begin{align*} 
&\int_{M_{\hat{t}}} \Phi \, \exp \Big ( (t_0-\hat{t})^{\frac{1}{2}}-\frac{e^{-K (t_0-\hat{t})} \, H}{200 \, H_1} \Big ) \, d\mu \\ 
&\leq \int_{M_{\tilde{t}}} \Phi \, \exp \Big ( (t_0-\tilde{t})^{\frac{1}{2}}-\frac{e^{-K (t_0-\tilde{t})} \, H}{200 \, H_1} \Big ) \, d\mu \\ 
&+ C \, |M_{\tilde{t}}| \, (t_0-\tilde{t}) + C(\Lambda) \, L^{-1} \, |M_{\tilde{t}}|, 
\end{align*}
provided that $\frac{1}{K} \geq t_0-\tilde{t} \geq t_0-\hat{t} \geq \frac{5}{9} \, H_1^{-2}$.
\end{theorem}

Note that the error term $C(\Lambda) \, L^{-1} \, |M_{\tilde{t}}|$ can be made arbitrarily small by choosing $L$ very large (depending on $\Lambda$). \\

\textbf{Proof.} 
By Proposition \ref{monotonicity.formula.in.manifolds}, the quantity 
\[\int_{M_t} \Phi \, \exp \Big ( (t_0-t)^{\frac{1}{2}}-\frac{e^{-K (t_0-t)} \, H}{200 \, H_1} \Big ) \, d\mu\] 
increases at a rate of at most $C \, |M_t|$ in between surgery times. Moreover, Proposition \ref{monotonicity.during.single.surgery.riemannian.case} implies that, during each surgery, the quantity 
\[\int_{M_t} \Phi \, \exp \Big ( (t_0-t)^{\frac{1}{2}}-\frac{e^{-K (t_0-t)} \, H}{200 \, H_1} \Big ) \, d\mu\] 
increases by at most $C(\Lambda) \, H_1^{-2}$, provided that $\frac{1}{K} \geq t_0-t \geq \frac{5}{9} \, H_1^{-2}$. On the other hand, the surface area decreases by at least $\frac{1}{10} \, L \, H_1^{-2}$ during each surgery. Hence, there are at most $C \, L^{-1} \, H_1^2 \, |M_{\tilde{t}}|$ surgeries after time $\tilde{t}$. Consequently, the quantity considered above increases by at most $C \, |M_{\tilde{t}}| \, (t_0-\tilde{t}) + C(\Lambda) \, L^{-1} \, |M_{\tilde{t}}|$ between time $\tilde{t}$ and time $\hat{t}$. \\

As a consequence of Theorem \ref{monotonicity.riemannian.case}, we obtain an analogue of Theorem \ref{regularity} in the Riemannian setting. 

\begin{theorem}
\label{regularity.riemannian.case}
Fix an open interval $I$ and a compact interval $J \subset I$. Moreover, suppose that $\bar{\mathcal{M}}$ is an embedded smooth solution of mean curvature flow in a Riemannian three-manifold which is defined for $t \in I$. Moreover, suppose that $\mathcal{M}^{(j)}$ is a sequence of mean curvature flows with surgery in the same Riemannian three-manifold, each of which is defined for $t \in I$. We assume that each surgery of the flow $\mathcal{M}^{(j)}$ involves performing $\Lambda$-surgery on an $(\hat{\alpha},\hat{\delta},\varepsilon,L)$-neck of size $r \in [\frac{1}{2H_1^{(j)}},\frac{1}{H_1^{(j)}}]$, where $H_1^{(j)} \to \infty$. We assume further that the surgery parameters satisfy $\hat{\delta} < \delta_*$, $\Lambda > \Lambda_*$, $\varepsilon < E_*(\Lambda)$, and that $L$ is chosen sufficiently large depending on $\Lambda$. Finally, we assume that the flows $\mathcal{M}^{(j)}$ converge to $\bar{\mathcal{M}}$ in the sense of geometric measure theory. Then, if $j$ is sufficiently large, the flow $\mathcal{M}^{(j)}$ is free of surgeries for all $t \in J$. Furthermore, the flows $\mathcal{M}^{(j)}$ converge smoothly to $\bar{\mathcal{M}}$ as $j \to \infty$. 
\end{theorem}

By combining Theorem \ref{regularity.riemannian.case} with results of Brian White \cite{White1},\cite{White2}, we are able to characterize the longtime behavior of mean curvature flow with surgery in Riemannian three-manifolds. This is discussed in \cite{Brendle-Huisken2}.


\begin{thebibliography}{99}
\bibitem{Brakke}
K.~Brakke, \textit{The motion of a surface by its mean curvature,} Princeton University Press (1978)

\bibitem{Brendle1}
S.~Brendle, \textit{A sharp bound for the inscribed radius under mean curvature flow,} to appear in Invent. Math.

\bibitem{Brendle2}
S.~Brendle, \textit{An inscribed radius estimate for mean curvature flow in Riemannian manifolds,} arxiv:1310.3439

\bibitem{Brendle-Huisken1}
S.~Brendle and G.~Huisken, \textit{Mean curvature flow with surgery of mean convex surfaces in $\mathbb{R}^3$,} to appear in Invent. Math.

\bibitem{Brendle-Huisken2}
S.~Brendle and G.~Huisken, \textit{Mean curvature flow with surgery of mean convex surfaces in three-manifolds,} preprint

\bibitem{Colding-Minicozzi}
T.~Colding and W.~Minicozzi, \textit{Generic mean curvature flow I: generic singularities,} Ann. of Math. 175, 755--833 (2012)

\bibitem{Ecker}
K.~Ecker, \textit{Regularity Theory for Mean Curvature Flow,} Birkh\"auser, Boston, 2004

\bibitem{Hamilton1}
R.~Hamilton, \textit{The formation of singularities in the Ricci flow,} Surveys in Differential Geometry, vol. II, 7--136, International Press, Somerville MA (1995)

\bibitem{Hamilton2} 
R.~Hamilton, \textit{Four-manifolds with positive isotropic curvature,} Comm. Anal. Geom. 5, 1--92 (1997)

\bibitem{Haslhofer-Kleiner1}
R.~Haslhofer and B.~Kleiner, \textit{Mean curvature flow of mean convex hypersurfaces,} arxiv:1304.0926

\bibitem{Haslhofer-Kleiner2}
R.~Haslhofer and B.~Kleiner, \textit{Mean curvature flow with surgery,} arxiv:1404.2332

\bibitem{Huisken}
G.~Huisken, \textit{Asymptotic behavior for singularities of the mean curvature flow,} J. Diff. Geom. 31, 285-299 (1990)

\bibitem{Huisken-Sinestrari}
G.~Huisken and C.~Sinestrari, \textit{Mean curvature flow with surgeries of two-convex hypersurfaces,} Invent. Math. 175, 137--221 (2009)

\bibitem{Perelman1} 
G.~Perelman, \textit{The entropy formula for the Ricci flow and its geometric applications,} arxiv:0211159

\bibitem{Perelman2}
G.~Perelman, \textit{Ricci flow with surgery on three-manifolds,} arxiv:0303109

\bibitem{White1}
B.~White, \textit{The size of the singular set in mean curvature flow of mean convex sets,} J. Amer. Math. Soc. 13, 665--695 (2000)

\bibitem{White2}
B.~White, \textit{The nature of singularities in mean curvature flow of mean convex sets,} J. Amer. Math. Soc. 16, 123--138 (2003)

\bibitem{White3}
B.~White, \textit{A local regularity theorem for mean curvature flow,} Ann. of Math. 161, 1487--1519 (2005)
\end{thebibliography}
\end{document}